\documentclass{amsart}

\usepackage{amsfonts,amssymb,verbatim,amsmath,amsthm,latexsym,textcomp,amscd}
\usepackage{latexsym,amsfonts,amssymb,epsfig,verbatim}
\usepackage{amsmath,amsthm,amssymb,latexsym,graphics,textcomp}
\usepackage{graphicx}
\usepackage{color}
\usepackage{url}

\input xy
\xyoption{all}

\begin{document}

\newtheorem{theorem}{Theorem}[section]
\newtheorem{prop}[theorem]{Proposition}
\newtheorem{lemma}[theorem]{Lemma}
\newtheorem{cor}[theorem]{Corollary}
\newtheorem{definition}[theorem]{Definition}
\newtheorem{conj}[theorem]{Conjecture}
\newtheorem{rmk}[theorem]{Remark}
\newtheorem{claim}[theorem]{Claim}
\newtheorem{defth}[theorem]{Definition-Theorem}

\newcommand{\boundary}{\partial}
\newcommand{\C}{{\mathbb C}}
\newcommand{\integers}{{\mathbb Z}}
\newcommand{\natls}{{\mathbb N}}
\newcommand{\ratls}{{\mathbb Q}}
\newcommand{\reals}{{\mathbb R}}
\newcommand{\proj}{{\mathbb P}}
\newcommand{\lhp}{{\mathbb L}}
\newcommand{\tube}{{\mathbb T}}
\newcommand{\cusp}{{\mathbb P}}
\newcommand\AAA{{\mathcal A}}
\newcommand\BB{{\mathcal B}}
\newcommand\CC{{\mathcal C}}
\newcommand\DD{{\mathcal D}}
\newcommand\EE{{\mathcal E}}
\newcommand\FF{{\mathcal F}}
\newcommand\GG{{\mathcal G}}
\newcommand\HH{{\mathcal H}}
\newcommand\II{{\mathcal I}}
\newcommand\JJ{{\mathcal J}}
\newcommand\KK{{\mathcal K}}
\newcommand\LL{{\mathcal L}}
\newcommand\MM{{\mathcal M}}
\newcommand\NN{{\mathcal N}}
\newcommand\OO{{\mathcal O}}
\newcommand\PP{{\mathcal P}}
\newcommand\QQ{{\mathcal Q}}
\newcommand\RR{{\mathcal R}}
\newcommand\SSS{{\mathcal S}}
\newcommand\TT{{\mathcal T}}
\newcommand\UU{{\mathcal U}}
\newcommand\VV{{\mathcal V}}
\newcommand\WW{{\mathcal W}}
\newcommand\XX{{\mathcal X}}
\newcommand\YY{{\mathcal Y}}
\newcommand\ZZ{{\mathcal Z}}
\newcommand\CH{{\CC\HH}}
\newcommand\PEY{{\PP\EE\YY}}
\newcommand\MF{{\MM\FF}}
\newcommand\RCT{{{\mathcal R}_{CT}}}
\newcommand\PMF{{\PP\kern-2pt\MM\FF}}
\newcommand\FL{{\FF\LL}}
\newcommand\PML{{\PP\kern-2pt\MM\LL}}
\newcommand\GL{{\GG\LL}}
\newcommand\Pol{{\mathcal P}}
\newcommand\half{{\textstyle{\frac12}}}
\newcommand\Half{{\frac12}}
\newcommand\Mod{\operatorname{Mod}}
\newcommand\Area{\operatorname{Area}}
\newcommand\ep{\epsilon}
\newcommand\hhat{\widehat}
\newcommand\Proj{{\mathbf P}}
\newcommand\U{{\mathbf U}}
 \newcommand\Hyp{{\mathbf H}}
\newcommand\D{{\mathbf D}}
\newcommand\Z{{\mathbb Z}}
\newcommand\R{{\mathbb R}}
\newcommand\Q{{\mathbb Q}}
\newcommand\E{{\mathbb E}}
\newcommand\til{\widetilde}
\newcommand\length{\operatorname{length}}
\newcommand\tr{\operatorname{tr}}
\newcommand\gesim{\succ}
\newcommand\lesim{\prec}
\newcommand\simle{\lesim}
\newcommand\simge{\gesim}
\newcommand{\simmult}{\asymp}
\newcommand{\simadd}{\mathrel{\overset{\text{\tiny $+$}}{\sim}}}
\newcommand{\ssm}{\setminus}
\newcommand{\diam}{\operatorname{diam}}
\newcommand{\pair}[1]{\langle #1\rangle}
\newcommand{\T}{{\mathbf T}}
\newcommand{\inj}{\operatorname{inj}}
\newcommand{\pleat}{\operatorname{\mathbf{pleat}}}
\newcommand{\short}{\operatorname{\mathbf{short}}}
\newcommand{\vertices}{\operatorname{vert}}
\newcommand{\collar}{\operatorname{\mathbf{collar}}}
\newcommand{\bcollar}{\operatorname{\overline{\mathbf{collar}}}}
\newcommand{\I}{{\mathbf I}}
\newcommand{\tprec}{\prec_t}
\newcommand{\fprec}{\prec_f}
\newcommand{\bprec}{\prec_b}
\newcommand{\pprec}{\prec_p}
\newcommand{\ppreceq}{\preceq_p}
\newcommand{\sprec}{\prec_s}
\newcommand{\cpreceq}{\preceq_c}
\newcommand{\cprec}{\prec_c}
\newcommand{\topprec}{\prec_{\rm top}}
\newcommand{\Topprec}{\prec_{\rm TOP}}
\newcommand{\fsub}{\mathrel{\scriptstyle\searrow}}
\newcommand{\bsub}{\mathrel{\scriptstyle\swarrow}}
\newcommand{\fsubd}{\mathrel{{\scriptstyle\searrow}\kern-1ex^d\kern0.5ex}}
\newcommand{\bsubd}{\mathrel{{\scriptstyle\swarrow}\kern-1.6ex^d\kern0.8ex}}
\newcommand{\fsubeq}{\mathrel{\raise-.7ex\hbox{$\overset{\searrow}{=}$}}}
\newcommand{\bsubeq}{\mathrel{\raise-.7ex\hbox{$\overset{\swarrow}{=}$}}}
\newcommand{\tw}{\operatorname{tw}}
\newcommand{\base}{\operatorname{base}}
\newcommand{\trans}{\operatorname{trans}}
\newcommand{\rest}{|_}
\newcommand{\bbar}{\overline}
\newcommand{\UML}{\operatorname{\UU\MM\LL}}
\newcommand{\EL}{\mathcal{EL}}
\newcommand{\tsum}{\sideset{}{'}\sum}
\newcommand{\tsh}[1]{\left\{\kern-.9ex\left\{#1\right\}\kern-.9ex\right\}}
\newcommand{\Tsh}[2]{\tsh{#2}_{#1}}
\newcommand{\qeq}{\mathrel{\approx}}
\newcommand{\Qeq}[1]{\mathrel{\approx_{#1}}}
\newcommand{\qle}{\lesssim}
\newcommand{\Qle}[1]{\mathrel{\lesssim_{#1}}}
\newcommand{\simp}{\operatorname{simp}}
\newcommand{\vsucc}{\operatorname{succ}}
\newcommand{\vpred}{\operatorname{pred}}
\newcommand\fhalf[1]{\overrightarrow {#1}}
\newcommand\bhalf[1]{\overleftarrow {#1}}
\newcommand\sleft{_{\text{left}}}
\newcommand\sright{_{\text{right}}}
\newcommand\sbtop{_{\text{top}}}
\newcommand\sbot{_{\text{bot}}}
\newcommand\sll{_{\mathbf l}}
\newcommand\srr{_{\mathbf r}}
\newcommand\geod{\operatorname{\mathbf g}}
\newcommand\mtorus[1]{\boundary U(#1)}
\newcommand\A{\mathbf A}
\newcommand\Aleft[1]{\A\sleft(#1)}
\newcommand\Aright[1]{\A\sright(#1)}
\newcommand\Atop[1]{\A\sbtop(#1)}
\newcommand\Abot[1]{\A\sbot(#1)}
\newcommand\boundvert{{\boundary_{||}}}
\newcommand\storus[1]{U(#1)}
\newcommand\Momega{\omega_M}
\newcommand\nomega{\omega_\nu}
\newcommand\twist{\operatorname{tw}}
\newcommand\modl{M_\nu}
\newcommand\MT{{\mathbb T}}
\newcommand\Teich{{\mathcal T}}
\renewcommand{\Re}{\operatorname{Re}}
\renewcommand{\Im}{\operatorname{Im}}

\title{Local and infinitesimal rigidity of simply connected negatively curved manifolds}

\author{Kingshook Biswas}
\address{RKM Vivekananda University, Belur Math, WB-711 202, India. Email: kingshook@rkmvu.ac.in}

\begin{abstract}
The boundary at infinity $\partial X$ of a CAT(-1) space $X$ carries a natural
family of metrics called visual metrics. These metrics are Moebius equivalent to each other (i.e. the
metric cross-ratios all coincide), and in particular conformal to each other, so the notions of
Moebius and conformal maps between boundaries of CAT(-1) spaces are well-defined, independent of choices
of visual metrics.

We consider a complete simply connected Riemannian manifold $(X, g_0)$ with sectional
curvatures bounded above by $-1$. The manifold $(X, g_0)$ is a CAT(-1) space and we denote its boundary
at infinity by $\partial_{g_0} X$. For a Riemannian metric $g$ on $X$ such that 
the tensor $g - g_0$ is compactly supported, and
such that all sectional curvatures of $g$ are bounded above by $-1$, the manifold $(X, g)$ is a
CAT(-1) space and the identity mapping $id : (X, g_0) \to (X, g)$ extends to a conformal homeomorphism
$\hat{id}_{g_0, g} : \partial_{g_0} X \to \partial_{g} X$.

Let $(g_t)_{0 \leq t \leq 1}$ be a $1$-parameter family of Riemannian metrics on $X$ depending
smoothly on $t$ such that the tensors $g_t - g_0, 0 \leq t \leq 1$, are supported in a fixed compact of $X$,
and such that all the metrics $g_t$ have sectional curvatures bounded above by $-1$. We show that if the
 all the boundary maps $\hat{id}_{g_0,g_t} : \partial_{g_0}X \to \partial_{g_t}X$
are Moebius, then the deformation is trivial, i.e. there is a $1$-parameter family of diffeomorphisms $f_t : X \to X$
such that $f^*_t g_t = g_0$.

We also show that given a compact $K \subset X$ and $0 < \alpha < 1$, there is an $\epsilon > 0$ such that if $g$ is a
Riemannian metric with sectional curvatures bounded above by $-1$ such that the support of $g - g_0$ is contained in $K$
and the $C^{2,\alpha}$ norm of $g - g_0$ is less than $\epsilon$, then if $\hat{id}_{g_0,g}$ is Moebius and $Vol_g(K) = Vol_{g_0}(K)$,
then $g$ is isometric to $g_0$.
\end{abstract}

\bigskip

\maketitle

\tableofcontents

\section{Introduction}

\medskip

The problems we consider in this article are motivated by rigidity
results for negatively curved manifolds. Two cases have been intensively studied, namely
closed negatively curved manifolds, and compact negatively curved manifolds with convex boundary. In both cases
rigidity results have been obtained, to the effect that metric deformations preserving some form of
a "length spectrum" are trivial, i.e. are isometric. For closed negatively curved manifolds, the role of length spectrum
is played by the {\it marked length spectrum}, i.e. the function on free homotopy classes of closed curves which assigns to
a homotopy class the length of the unique closed geodesic in that class, while for compact negatively curved manifolds
with boundary, the role of length spectrum is played by the {\it boundary distance function}, i.e. the function which
assigns to pairs of points on the boundary the geodesic distance between them.

\medskip

 To put the corresponding results in
a general context, we may consider the moduli space $M(X)$ of negatively curved metrics on a manifold $X$, i.e. the quotient of
the space of negatively curved metrics on $X$ by the natural action of the group of diffeomorphisms of $X$. We are then given
a "length spectrum map" $L : M(X) \to C^{\mathbb{R}^+}$ which assigns to a negatively curved metric $g$ a length function $L(g) : C \to \mathbb{R}^+$,
where $C$ is a space parametrizing a certain collection of geodesics (in the two cases mentioned above, $C$ would be the set of free
homotopy classes of closed curves and the set of pairs of distinct points on the boundary respectively). There are then three rigidity
questions one may pose: global rigidity (injectivity of the map $L$), 
local rigidity (local injectivity of $L$) and
infinitesimal rigidity (injectivity of the differential of $L$).

\medskip


We mention briefly some results obtained in the two cases previously mentioned. Guillemin and Kazhdan \cite{guilleminkazhdan}
proved an infinitesimal rigidity result for closed negatively curved surfaces. Otal proved that global rigidity holds for
closed negatively curved surfaces \cite{otal2} (giving an affirmative answer to the "marked length spectrum rigidity"
problem of Burns-Katok \cite{burnskatok} in dimension $2$), and also for compact negatively curved surfaces with convex boundary \cite{otal3}.
Croke and Sharafutdinov have proven that infinitesimal rigidity holds for closed negatively curved $n$-manifolds
\cite{crokesharafutdinov}, and Croke-Dairbekov-Sharafutdinov \cite{crokedairbekovsharafutdinov} 
have also proven a local rigidity result for compact negatively curved $n$-manifolds with convex boundary.

\medskip

We prove local and infinitesimal rigidity results
in a third case not previously considered, namely that of simply connected complete negatively curved
manifolds. This is similar to the second case, only now the boundary is at infinity. In this case given a complete
negativey curved metric $g_0$ on a simply connected manifold $X$, say with sectional curvatures bounded above by $-1$ so
that $X$ with the distance function induced by $g_0$ is a CAT(-1) space, while there is no natural notion of a length spectrum
for the metric $g_0$, there is, given a compactly supported deformation $g$ of $g_0$,
a well-defined notion of a {\it relative length spectrum} for the pair $(g_0, g_1)$. This is given by a
function $S_{g_0}(g) : \partial^2_{g_0} X \to \mathbb{R}$ of pairs of distinct points on the boundary at infinity
$\partial_{g_0}X$ of $X$ with respect to $g_0$,
called the {\it integrated Schwarzian} of $g$ with respect to $g_0$, which is a renormalized version
of the boundary distance function, measuring the difference between the $g_0$ and $g$ distances between
pairs of points on the boundary at infinity.

\medskip

The integrated Schwarzian was first introduced in \cite{biswas3} in the
study of Moebius and conformal maps between boundaries of CAT(-1) spaces. We recall that the boundary at
infinity $\partial X$ of a CAT(-1) space $X$ carries a natural
family of metrics called visual metrics. These metrics are Moebius equivalent to each other (i.e. the
metric cross-ratios all coincide), and in particular conformal to each other, hence the notions of
Moebius and conformal maps between boundaries of CAT(-1) spaces are well-defined, independent of choices
of visual metrics (for basic properties of CAT(-1) spaces and their boundaries see \cite{bourdon1}, \cite{bourdon2}).
Given a conformal map between two boundaries $f : \partial X \to \partial Y$, the integrated Schwarzian of $f$ is
a function $S(f) : \partial^2 X \to \mathbb{R}$ which measures the deviation of the conformal map $f$ from being Moebius. In the
case of a pair of negatively curved metrics $g_0, g_1$ as above with $g - g_0$ compactly supported, the identity map $id : (X, g_0) \to (X, g)$
extends in fact to a conformal homeomorphism $\hat{id}_{g_0,g} : \partial_{g_0}X \to \partial_g X$. The integrated Schwarzian of $g$ with respect to
$g_0$ is then defined to be the integrated Schwarzian of the conformal map $\hat{id}_{g_0,g}$. It turns out that the
integrated Schwarzian vanishes if the boundary
map $\hat{id}_{g_0,g} : \partial_{g_0} X \to \partial_{g}X$ is Moebius, and conversely in the presence
of a lower curvature bound for $g_0$, $\hat{id}_{g_0,g}$ is Moebius if $S_{g_0}(g)$ vanishes.

\medskip

We prove the following infinitesimal rigidity result:

\medskip

\begin{theorem} \label{isospectral} Let $(X, g_0)$ be a complete simply connected Riemannian manifold with
sectional curvatures bounded above by $-1$. Let $(g_t)_{0 \leq t \leq 1}$ be a $1$-parameter family of Riemannian
metrics on $X$ such that:

\smallskip

\noindent (1) The symmetric $(0,2)$-tensors $g_t - g_0, 0 \leq t \leq 1$ are compactly supported with supports contained in
 a fixed compact $C \subset X$.

\smallskip

\noindent (2) The sectional curvatures of the metrics $g_t, 0 \leq t \leq 1$ are bounded above by $-1$.

\smallskip

\noindent (3) The metrics $g_t, 0 \leq t \leq 1$ depend smoothly on the parameter $t$, i.e. the map $[0,1] \times X \to T^*X^{\odot 2}, (t,x) \mapsto g_t(x)$
is smooth.

\smallskip

\noindent (4) All the boundary maps $\hat{id}_{g_0,g_t} : \partial_{g_0}X \to \partial_{g_t}X, 0 \leq t \leq 1$, are Moebius.

\medskip

Then there is a $1$-parameter family of diffeomorphisms $f_t : X \to X$ such that $f^*_t g_t = g_0$. Moreover there is a
compact $K \subset X$ such that $f_t = id_X$ on $X - K$ for $0 \leq t \leq 1$.
\end{theorem}

\medskip

There are two main steps in the proof of the above theorem. The first is to prove a formula for the first variation
of the integrated Schwarzian,

$$
\frac{d}{dt} 2S_{g_0}(g_t) = I_{g_t}(\dot{g_t}) \circ \hat{id}_{g_0, g_t}
$$

\noindent where $I_{g_t}$ denotes the {\it ray transform} of the metric $g_t$ and $\dot{g_t}$ the symmetric $(0,2)$-tensor $\frac{d}{dt}(g_t)$.
Here by the ray transform of a metric $g$ we mean the map which assigns to a compactly supported symmetric $(0,2)$-tensor $u$
the function $I_g(u) : \partial^2 X \to \mathbb{R}$ on the space of bi-infinite geodesics of $X$ obtained by integrating
$u$ along geodesics (when $u$ is thought of as a function on the unit tangent bundle). This
formula is an analogue of a formula for the first variation of the boundary distance function in the case of a compact manifold
with boundary, and indeed is proved by passing to the limit in this formula.

\medskip

The second step is to prove that the kernel of the ray transform consists of {\it exact} symmetric $(0,2)$-tensors, i.e.
symmetric $(0,2)$-tensors arising as Lie derivatives of the metric with respect to vector fields.
The analogous statement in the case of a compact manifold
with boundary is known, and again the statement in the simply connected case is proved by passing to the limit.
Combining these two steps, it follows that each $\dot{g_t}$ is
exact, and integrating the corresponding vector fields gives the required family of diffeomorphisms $(f_t)_{0 \leq t \leq 1}$.

\medskip

We also prove a local rigidity result:

\begin{theorem} \label{localrigidity} Let $(X, g_0)$ be a complete simply
connected Riemannian manifold with sectional curvature
bounded above by $-1$. Given a compact $K \subset X$ and
$0 < \alpha < 1$, there is an $\epsilon > 0$ such that the following
holds:

Let $g$ be a metric with sectional curvatures bounded
above by $-1$ such that the support of $g - g_0$ is contained in $K$, and such that
the $C^{2,\alpha}$ norm of $g - g_0$ is less than $\epsilon$.
If the boundary map $\hat{id}_{g_0,g} : \partial_{g_0}X \to \partial_g X$ is
Moebius and $Vol_{g_0}(K) = Vol_g(K)$ then $g$ is isometric to $g_0$.
\end{theorem}

\medskip

The proof of the above theorem follows along lines very similar to Croke-Dairbekov-Sharafutdinov's proof of local rigidity for
compact manifolds with convex boundary, and is essentially an adaptation of their proof to the case at hand.

\medskip

The paper is organized as follows. In section 2 we recall background material on Moebius maps, conformal maps
and the integrated Schwarzian, in the context of general CAT(-1) spaces. In section 3 we 
consider the integrated Schwarzian in the case of a compactly supported 
deformation of a complete simply connected negatively curved Riemannian manifold,
and derive the variational formula described above. In section 4 we prove the assertion 
about the kernel of the ray transform.
Finally in section 5 we put these ingredients together to prove Theorems \ref{isospectral} 
and \ref{localrigidity}.

\medskip

\section{Moebius maps, conformal maps and the integrated Schwarzian}

\medskip

The material in this section is taken from \cite{biswas3}.

\medskip

\begin{definition} A homeomorphism between metric spaces $f :
(Z_1, \rho_1) \to (Z_2, \rho_2)$ with no isolated points is said to be {\it conformal} if
for all $\xi \in Z_1$, the limit
$$
df_{\rho_1, \rho_2}(\xi) := \lim_{\eta \to \xi} \frac{\rho_2(f(\xi),
f(\eta))}{\rho_1(\xi, \eta)}
$$
exists and is positive. The positive function $df_{\rho_1,
\rho_2}$ is called the derivative of $f$ with respect to $\rho_1, \rho_2$.
We say $f$ is {\it $C^1$ conformal} if its derivative is continuous.

\medskip

Two metrics $\rho_1, \rho_2$ inducing the same topology on a set
$Z$, such that $Z$ has no isolated points,
are said to be conformal (respectively $C^1$ conformal) if the
map $id_Z : (Z, \rho_1) \to (Z, \rho_2)$ is conformal
(respectively $C^1$ conformal). In this case we denote the
derivative of the identity map by $\frac{d\rho_2}{d\rho_1}$.
\end{definition}

\medskip

\begin{definition} Let $Z$ be a set with at least four points. For a metric $\rho$ on
$Z$ we define the metric cross-ratio with respect to $\rho$ of a quadruple of distinct
points $(\xi, \xi', \eta, \eta')$ of $Z$ by
$$
[\xi \xi' \eta \eta']_{\rho} := \frac{\rho(\xi, \eta) \rho(\xi', \eta')}{\rho(\xi,
\eta')\rho(\xi', \eta)}
$$
\end{definition}

\medskip

\begin{definition} A map between metric spaces $f : (Z_1, \rho_1) \to (Z_2, \rho_2)$ is said to be Moebius if
it preserves metric cross-ratios. A map $f$ is locally Moebius if every $\xi \in Z_1$ has a neighbourhood $U$ such that
$f_{|U}$ is Moebius. Two metrics $\rho_1, \rho_2$ on a set $Z$ are Moebius equivalent if the identity
map $id : (Z, \rho_1) \to (Z, \rho_2)$ is Moebius.
\end{definition}

\medskip

We recall the following facts:

\begin{prop} We have the following:

\smallskip

\noindent (1) A locally Moebius map between metric spaces with no isolated points is $C^1$ conformal.

\smallskip

\noindent (2) A Moebius map $f : (Z_1, \rho_1) \to (Z_2, \rho_2)$ between metric spaces with
no isolated points satisfies the "geometric mean-value theorem":

$$
\rho_2(f(\xi), f(\eta))^2 = df_{\rho_1,\rho_2}(\xi) df_{\rho_1,\rho_2}(\eta) \rho_1(\xi, \eta)^2
$$

\end{prop}

\medskip

Let $X$ be a CAT(-1) space. We recall that the boundary at infinity of $X$, denoted by $\partial X$, comes equipped with
a family of metrics $(\rho_x)_{x \in X}$ called visual metrics, defined by $\rho_x(\xi, \eta) = \exp(-(\xi|\eta)_x)$, where
$(.|.)_x$ denotes the Gromov inner product (extended to the boundary) with respect to the basepoint $x$. These metrics are all Moebius
equivalent, so the notions of Moebius and conformal maps between boundaries of CAT(-1) spaces $f : \partial X \to \partial Y$ are
well-defined independent of the choices of visual metrics on $\partial X$ and $\partial Y$. For $(\xi, \eta) \in \partial^2 X$,
we denote the bi-infinite geodesic in $X$ with endpoints $\xi, \eta$ also by $(\xi, \eta)$.

\medskip

\begin{definition} Let $f : \partial X \to \partial Y$ be a
conformal map between boundaries of CAT(-1) spaces equipped with
visual metrics. The integrated Schwarzian of $f$ is the function
$S(f) : \partial^2 X \to \mathbb{R}$ defined by
$$
S(f)(\xi, \eta) := \log (df_{\rho_x, \rho_y}(\xi) df_{\rho_x,
\rho_y}(\eta)) \ \ (\xi,\eta) \in \partial^2 X
$$
where $x,y$ are any two points $x \in (\xi, \eta), y \in (f(\xi),
f(\eta))$ (it is easy to see that the quantity defined above is
independent of the choices of $x$ and $y$).
\end{definition}

\medskip

It is easy to see (using the geometric mean-value theorem) that the integrated Schwarzian of a
Moebius map is identically zero. We recall that a CAT(-1) space is said to be {\it geodesically complete}
if every geodesic segment can be extended to a bi-infinite geodesic. We have the following from \cite{biswas3}:

\medskip

\begin{theorem} \label{moebiusconformal} Let $X$ be a simply
connected complete Riemannian manifold with sectional curvatures
satisfying $-b^2 \leq K \leq -1$ for some $b \geq 1$, and let $Y$ be a proper geodesically complete CAT(-1) space.
A $C^1$ conformal map $f : U \subset \partial X \to V \subset \partial
Y$ between open subsets of $\partial X, \partial Y$ is Moebius on $U$ (i.e. preserves cross-ratios)
if and only if its integrated Schwarzian $S(f)$ vanishes on $\partial^2
U$. Two $C^1$ conformal maps $f, g : U \subset \partial X \to V \subset \partial
Y$ differ by post-composition with a Moebius map $h : V \to V$ if
and only if their integrated Schwarzians coincide, $S(f) = S(g)$.
\end{theorem}

\medskip

We recall that the Busemann function of a CAT(-1) space $X$ is the function $B : \partial X \times X \times X \to \mathbb{R}$
defined by
$$
B(\xi, x, y) := \lim_{a \to \xi} \left(d_X(x,a) - d_X(y,a)\right)
$$
where $a \in X$ converges radially towards $\xi$.

\medskip

We have the following formula for the derivatives of visual metrics on $\partial X$:

$$
\frac{d\rho_y}{d\rho_x}(\xi) = \exp(B(\xi,x,y)) \ \xi \in \partial X, x,y \in X.
$$

\medskip

\section{The integrated Schwarzian for compactly supported deformations}

\medskip

Let $(X, g_0)$ be a complete simply connected Riemannian manifold with sectional curvatures
bounded above by $-1$. Let $g$ be a Riemannian metric on $X$ with sectional curvatures bounded above by $-1$
such that the symmetric $(0,2)$-tensor $g - g_0$ is compactly supported. Note that $(X, g_0)$ and $(X, g)$ are both
CAT(-1) spaces. We denote the corresponding distance functions on $X$ by $d_{g_0}, d_g$, the
boundaries by $\partial_{g_0}X, \partial_g X$ and the visual metrics by $\rho_{x, g_0}, \rho_{x, g}, x \in X$.

\medskip

\begin{lemma} The identity map $id : (X, g_0) \to (X, g)$ extends to a locally Moebius homeomorphism
$\hat{id}_{g_0,g} : \partial_{g_0} X \to \partial_{g} X$ (in particular the boundary map is conformal).
\end{lemma}

\medskip

\noindent{\bf Proof:} Clearly $id : (X, g_0) \to (X, g)$ is bi-Lipschitz, hence extends to a
homeomorphism $\hat{id}_{g_0,g} : \partial_{g_0} X \to \partial_{g} X$. Fix a basepoint $x_0 \in X$ and let
$B \subset X$ be a large $g_0$-ball around $x_0$ containing the support of $g - g_0$.
Given a point $\xi_0 \in \partial_{g_0}X$, we may choose a neighbourhood $U$ of $\xi_0$ small enough
such that for all $\xi, \eta \in U, \xi \neq \eta$, the bi-infinite $g_0$-geodesic with
endpoints $\xi, \eta$ lies outside $B$. We may also choose $x_1$ lying along the $g_0$-geodesic ray joining
$x_0$ to $\xi_0$ far enough from $x_0$ such that for all $\xi \in U$, the $g_0$-geodesic ray joining $x_1$ to
$\xi$ lies outside $B$. Since $g_0$-geodesics lying outside $B$ are also $g_1$-geodesics, it follows that
$\hat{id}_{g_0, g} : (U, \rho_{x_1, g_0}) \to (V = \hat{id}_{g_0, g}(U), \rho_{x_1, g})$ is an isometry. Since the
metrics $\rho_{x_0, g_0}, \rho_{x_1, g_0}$ are Moebius equivalent, as are the metrics $\rho_{x_0, g}, \rho_{x_1, g}$,
it follows that $\hat{id}_{g_0, g} : (U, \rho_{x_0, g_0}) \to (V, \rho_{x_0,g})$ is Moebius. $\diamond$

\medskip

\begin{definition} The integrated Schwarzian of $g$ with respect to $g_0$ is defined to be the integrated
Schwarzian of the above conformal boundary map, $S_{g_0}(g) := S(\hat{id}_{g_0,g}) : \partial^2_{g_0} X \to \mathbb{R}$.
\end{definition}

\medskip

The following lemma says that the integrated Schwarzian may be viewed as a renormalized limit of boundary distance functions:

\medskip

\begin{lemma} \label{schwarzlimit} For $(\xi, \eta) \in \partial^2_{g_0} X$, the following limit exists and equals the integrated Schwarzian:
$$
\lim_{p,q\in (\xi, \eta), p \to \xi, q\to \eta} (d_{g}(p,q) - d_{g_0}(p,q)) = S_{g_0}(g)(\xi,\eta)
$$
\end{lemma}

\medskip

\noindent{\bf Proof:} Let $f = \hat{id}_{g_0,g}$. Let $\gamma_0 : \mathbb{R} \to X$ be the unique bi-infinite $g_0$-geodesic with endpoints
$\gamma_0(-\infty) = \eta, \gamma_0(+\infty) = \xi$ in $\partial_{g_0} X$, and let $\gamma : \mathbb{R} \to X$ be the
unique bi-infinite $g$-geodesic with endpoints $\gamma(-\infty) = \eta' := f(\eta),
\gamma(+\infty) = \xi' := f(\xi)$ in $\partial_g X$. By definition, $d_{g_0}(\gamma_0(t), \gamma(t))$ is bounded as
$t \to -\infty$ and as $t \to +\infty$. Since $g - g_0$ is compactly supported, for $R > 0$ large enough, $\gamma^+_R := \gamma_{|[R,+\infty)}$
and $\gamma^-_R := \gamma_{|(-\infty, -R]}$ are $g_0$-geodesic rays with endpoints $\xi, \eta$ respectively in $\partial_{g_0} X$.

\medskip

Let $x_t = \gamma_0(t), y_t = \gamma_0(-t)$ for $t > 0$. For $t > 0$ large enough, let $x'_t \in \gamma^+_R, y'_t \in \gamma^-_R$ be such that
$B(\xi, x_t, x'_t) = 0, B(\eta, y_t, y'_t) = 0$, where $B : \partial_{g_0} X \times X \times X \to \mathbb{R}$ is the Busemann function
of $(X, g_0)$. Note that $d_{g_0}(x_t, x'_t) \to 0, d_{g_0}(y_t, y'_t) \to 0$ as $t \to +\infty$ by exponential convergence of asymptotic
geodesic rays in the CAT(-1) space $(X, g_0)$.

\medskip

From the formula for the derivatives of visual metrics, we have
$$
\frac{d\rho_{x_t,g_0}}{d\rho_{x'_t,g_0}} = \frac{d\rho_{y_t,g_0}}{d\rho_{y'_t,g_0}} = 1
$$
Since in the previous Lemma we saw that the restrictions of $f$ to small neighbourhoods $U, U'
$ of $\xi, \eta$ respectively were isometries $f : (U, \rho_{x'_t, g_0}) \to (f(U), \rho_{x'_t, g}),
f : (U, \rho_{y'_t, g_0}) \to (f(U), \rho_{y'_t, g})$ for $t$ large enough, it follows that
$$
df_{\rho_{x_t, g_0}, \rho_{x'_t,g}}(\xi) = df_{\rho_{y_t, g_0}, \rho_{y'_t,g}}(\eta) = 1
$$
Hence by definition of the integrated Schwarzian and the chain rule we have
\begin{align*}
S_{g_0}(g)(\xi,\eta) & = \log df_{\rho_{x_t, g_0}, \rho_{x'_t,g}}(\xi) + \log df_{\rho_{x_t,g_0}, \rho_{x'_t,g}}(\eta) \\
                     & = \log df_{\rho_{x_t,g_0}, \rho_{x'_t,g}}(\eta) \\
                     & = \log \frac{d\rho_{y'_t,g}}{d\rho_{x'_t,g}}(\eta) + \log \frac{d\rho_{x_t, g_0}}{d\rho_{y_t,g_0}}(\eta) \\
                     & = d_{g}(x'_t, y'_t) - d_{g_0}(x_t,y_t) \\
                     & = d_{g}(x_t,y_t) - d_{g_0}(x_t,y_t) + o(1) \\
\end{align*}

\noindent as $t \to +\infty$, since for $t$ large, we have $d_{g}(x'_t, x_t) = d_{g_0}(x'_t,x_t) \to 0, d_{g}(y'_t, y_t) = d_{g_0}(y'_t,y_t) \to 0$.
The result follows. $\diamond$

\medskip

We also have the following elementary lemma:

\medskip

\begin{lemma}\label{distbdd} There is a constant $C > 0$ such that for any $p,q \in X$ we have
$$
|d_g(p,q) - d_{g_0}(p,q)| \leq C
$$
\end{lemma}

\medskip

\noindent{\bf Proof:} Let $B$ be an open $g_0$-ball containing the
support of $g - g_0$. Since $\overline{B}$
is compact we may assume
that $p,q$ are not both contained in $\overline{B}$. We
give the proof only in the case when both
$p$ and $q$ lie outside $B$, the argument being similar in
the other case when one of the two points
lies outside $B$ and the other inside $B$.

\medskip

Let $\gamma_0, \gamma$ be $g_0$ and $g$ geodesics respectively
joining $p,q$.
If either one of the curves $\gamma_0, \gamma$ does not intersect $B$
then neither does the other (since a $g_0$-geodesic lying entirely outside $B$ is
also a $g$-geodesic and vice-versa, and $g_0, g$-geodesics joining points are unique)
and $d_g(p,q) = d_{g_0}(p,q)$. Otherwise choose points $a,b$ on $\gamma_0 \cap \partial B$ and points
$a',b'$ on $\gamma \cap \partial B$ such that the $g_0$-geodesic segments $[p,a]$ and $[b,q]$
and the $g$-geodesic segments
$[p,a'], [b',q]$ are all disjoint from $B$, then
\begin{align*}
|d_{g}(p,a') - d_{g_0}(p,a)| & = |d_{g_0}(p,a') - d_{g_0}(p,a)| \leq d_{g_0}(a',a), \\
|d_{g}(b',q) - d_{g_0}(b,q)| & = |d_{g_0}(b',q) - d_{g_0}(b,q)| \leq d_{g_0}(b',b), \\
\end{align*}
while $|d_g(a',b') - d_{g_0}(a,b)|$ is bounded by $diam_g(B) + diam_{g_0}(B)$,
so that $|d_g(p,q) - d_{g_0}(p,q)| \leq diam_g(B) + 3 diam_{g_0}(B)$. $\diamond$

\medskip

We now consider a $1$-parameter family of Riemannian metrics $(g_t)_{0 \leq t \leq 1}$ on $X$
satisfying the following hypotheses:

\medskip

\noindent (1) There is a fixed compact $K$ containing the supports of the symmetric $(0,2)$-tensors $g_t - g_0, 0 \leq t \leq 1$.

\medskip

\noindent (2) The sectional curvatures of the metrics $g_t, 0 \leq t \leq 1$, are bounded above by $-1$.

\medskip

\noindent (3) The metrics $g_t, 0 \leq t \leq 1$ depend smoothly on the parameter $t$, i.e. the map $[0,1] \times X \to T^*X^{\odot 2}, (t,x) \mapsto g_t(x)$
is smooth.

\medskip

\begin{lemma} Fix $p,q \in X, p \neq q$. Let $a = d_{g_0}(p,q)$, and for $0 \leq t \leq 1$, let $\gamma_t : [0,a] \to X$ be the
unique $g_t$-geodesic segment with endpoints $p,q$. Then the map $[0,1] \times [0,a] \to X, (t,s) \mapsto \gamma_t(s)$ is smooth.
\end{lemma}

\medskip

\noindent{\bf Proof:} Let $\exp_t : T_p X \to X$ denote the exponential mapping of the metric $g_t$ based at $p$. By smooth
dependence of solutions to ODE's on initial conditions and on coefficients, the map $\Phi : [0,1] \times T_p X \to X, (t,v) \mapsto \exp_t(v)$
is smooth. Since the metrics $g_t$ are nonpositively curved, for each $t$ there is a unique $v_t \in T_p X$ such that $\Phi(t, v_t) = q$. Moreover
each map $\exp_t$ is a diffeomorphism, hence applying the Implicit Function Theorem to $\Phi$ it follows that the map $t \mapsto v_t$ is smooth.
Since $\gamma_t(s) = \Phi(t, (s/a)v_t)$, the lemma follows. $\diamond$

\medskip

For notational convenience, given a symmetric $(0,m)$-tensor field $u$ on $X$ and a tangent vector $\xi \in TX$, we denote
$u(\xi,\dots,\xi)$ by simply $u(\xi)$.

\medskip

\begin{lemma} \label{distvarn} With the same notation as above, we have:
$$
\frac{d}{dt}\left( \frac{d_{g_t}(p,q)^2 - d_{g_0}(p,q)^2}{d_{g_0}(p,q)}\right) = \int_{0}^{a} \dot{g_t}(\dot{\gamma_t}(s)) \ ds
$$
where $\dot{g_t}$ is the symmetric $(0,2)$-tensor $\frac{d}{dt}(g_t)$.
\end{lemma}

\medskip

\noindent{\bf Proof:} This formula may be found in \cite{sharafutdinov1}. We reproduce the proof for the benefit of the reader. It suffices
to prove the formula for $t = 0$. Let $\gamma_t(s) = (\gamma^i_t(s)), g_t = (g_{t,ij}), \dot{g_0} = (f_{ij})$ in local coordinates. We have (using
Einstein summation convention)
$$
\frac{d_{g_t}(p,q)^2}{d_{g_0}(p,q)} = \int_{0}^{a} g_{t,ij}(\gamma_t(s))\dot{\gamma}^i_t(s)\dot{\gamma}^j_t(s) \ ds
$$
Differentiating the above equality with respect to $t$ and putting $t = 0$ gives
$$
\frac{d}{dt}_{|t = 0}\left( \frac{d_{g_t}(p,q)^2 - d_{g_0}(p,q)^2}{d_{g_0}(p,q)}\right) =
\int_{0}^{a} f_{ij}(\gamma_0(s))\dot{\gamma}^i_0(s)\dot{\gamma}^j_0(s) \ ds + A
$$
where $A$ is the integral
$$
A = \int_{0}^{a} \left( \frac{\partial}{\partial x^k}g_{0,ij}(\gamma_0(s))\dot{\gamma}^i_0(s)\dot{\gamma}^j_0(s)\frac{\partial}{\partial t}_{|t = 0}\gamma^k_t(s)
+ 2g_{0,ij}(\gamma_0(s))\dot{\gamma}^i_0(s)\frac{\partial}{\partial t}_{|t=0}\dot{\gamma}^j_t(s) \right) ds
$$
Now the integral $A$ is equal to zero, since $\frac{\partial}{\partial t}\gamma_t(0) = \frac{\partial}{\partial t}\gamma_t(a) = 0$,
and the curve $\gamma_0$ is an extremal of the energy functional
$$
E_0(\gamma) = \int_{0}^{a} g_{0,ij}(\gamma(s))\dot{\gamma}^i(s)\dot{\gamma}^j(s) \ ds
$$
The lemma follows. $\diamond$

\medskip

Denote by $C^{\infty}_c(T^*X^{\odot 2})$ the space of smooth compactly supported 
symmetric $(0,2)$-tensors on $X$.

\medskip

\begin{definition} Let $g$ be a complete Riemannian metric on $X$ with sectional curvatures bounded above by $-1$. The
ray transform of $g$ is the linear map
\begin{align*}
I_g : C^{\infty}_{c}(T^*X^{\odot 2}) & \to C_c(\partial^2_g X) \\
                         f            & \mapsto \left( I_g(f) : (\xi, \eta) \mapsto \int_{-\infty}^{\infty} f(\dot{\gamma}_{(\xi, \eta)}(s)) \ ds\right) \\
\end{align*}
\noindent where $\gamma_{(\xi, \eta)}$ is the unique (up to translation) bi-infinite $g$-geodesic with unit $g$-speed and
endpoints $\gamma_{(\xi,\eta)}(-\infty) = \eta,
\gamma_{(\xi,\eta)}(+\infty) = \xi$.
The domain of $I_g$ is the space of smooth symmetric $(0,2)$-tensors on $X$ with compact support, and the range the space of
continuous functions on $\partial^2_g X$ with compact support.
\end{definition}

\medskip

We can now prove the first variation formula for the integrated Schwarzian:

\medskip

\begin{theorem} \label{schwarzvarn} With the same notation as above, for $(\xi, \eta) \in \partial^2_{g_0}X$ we have:
$$
\frac{d}{dt}2S_{g_0}(g_t)(\xi,\eta) = \left(I_{g_t}(\dot{g_t}) \circ \hat{id}_{g_0,g_t}\right) (\xi, \eta)
$$
\end{theorem}

\medskip

\noindent{\bf Proof:} For $R > 0$, let $p_R = \gamma_{(\xi, \eta)}(R), q_R = \gamma_{(\xi, \eta)}(-R)$ where
$\gamma_{(\xi, \eta)}$ is a bi-infinite unit speed $g_0$-geodesic with endpoints $\xi,\eta$. Let $\gamma_{t,R} : [-R,R] \to X$ be
the unique $g_t$-geodesic segment with endpoints $p_R, q_R$. Define real-valued functions on $[0,1]$ by
$$
h_R : t \mapsto \frac{d_{g_t}(p_R,q_R)^2 - d_{g_0}(p_R,q_R)^2}{d_{g_0}(p_R,q_R)}
$$
Then it is easy to see using Lemma \ref{schwarzlimit} and Lemma \ref{distbdd} that the pointwise limit as $R \to +\infty$ of the functions $h_R$
is the function $h : t \mapsto 2 S_{g_0}(g_t)(\xi,\eta)$. Moreover from Lemma \ref{distvarn} it follows that each $h_R$ is
differentiable with derivative
$$
h'_R : t \mapsto \int_{-R}^{R} \dot{g_t}(\dot{\gamma}_{t,R}(s)) \ ds
$$
By Lemma \ref{distvarn}, $\gamma_{t,R}$ depends smoothly on $t$, hence $h'_R$ is continuous and $h_R$ is in fact $C^1$.
Moreover as $R \to +\infty$, for each fixed $t$ it follows from a standard argument for CAT(-1) spaces that the
$g_t$-geodesic segments $\gamma_{t,R}$ converge uniformly on compacts to a bi-infinite
$g_t$-geodesic $\gamma_t : \mathbb{R} \to X$
with endpoints $\hat{id}_{g_0,g_t}(\xi), \hat{id}_{g_0,g_t}(\eta)$, which is unit speed (because $d_{g_t}(p_R, q_R)/d_{g_0}(p_R, q_R) \to 1$
as $R \to +\infty$). The same argument for CAT(-1) spaces also gives that
the convergence of $\gamma_{t,R}$ to $\gamma_t$ is uniform in $t$ (the upper bound on the
distance between $\gamma_{t,R}$ and $\gamma_t$ only depends on the upper sectional curvature bound for $g_t$, which is $-1$
independent of $t$, and on the visual distance $\rho_{x,g_t}(\hat{id}_{g_0, g_t}(\xi), \hat{id}_{g_0, g_t}(\eta))$, which is
bounded below by a positive constant independent of $t$ for a fixed basepoint $x \in X$). Again due to negative curvature,
$C^0$-convergence of $\gamma_{t,R}$ to $\gamma_t$ actually implies $C^1$-convergence of $\gamma_{t,R}$ to $\gamma_t$.
It follows that as
$R \to +\infty$, the functions $h'_R$ converge uniformly to the function
$$
f : t \mapsto \int_{-\infty}^{\infty} \dot{g_t}(\dot{\gamma}_t(s)) \ ds = \left(I_{g_t}(\dot{g_t}) \circ \hat{id}_{g_0,g_t}\right) (\xi, \eta)
$$
It follows that $h$ is $C^1$ with derivative equal to $f$. $\diamond$

\medskip

We will need the following lemma in the proof of the local rigidity result:

\medskip

\begin{lemma}\label{raybiggerschwarz} Let $g_1$ be a negatively curved metric on $X$ with $g_1 - g_0$ compactly supported,
then for any $(\xi, \eta) \in \partial^2_{g_0}X$ we have:
$$
I_{g_0}(g_1 - g_0) (\xi, \eta) \geq 2S_{g_0}(g_1)(\xi,\eta)
$$
\end{lemma}

\medskip

\noindent{\bf Proof:} For $p,q \in X$ and $i = 0,1$ let $\gamma^i_{p,q} : [0, d_{g_0}(p,q)] \to X$ be the unique
$g_i$-geodesic joining $p$ to $q$. Letting $p,q$ tend to $\xi,\eta$ radially along $g_0$ geodesics, we have:

\begin{align*}
I_{g_0}(g_1 - g_0) (\xi, \eta) & = \lim_{p \to \xi, q \to \eta} \left(\int_{0}^{d_{g_0}(p,q)} g_1((\gamma^{0})'(t)) dt
                                 - \int_{0}^{d_{g_0}(p,q)} g_0((\gamma^{0})'(t)) dt \right) \\
                               & \geq \lim_{p \to \xi, q \to \eta} \left(\int_{0}^{d_{g_0}(p,q)} g_1((\gamma^{1})'(t)) dt
                                 - \int_{0}^{d_{g_0}(p,q)} g_0((\gamma^{0})'(t)) dt \right) \\
                               & = \lim_{p \to \xi, q \to \eta} \left( \frac{d^2_{g_1}(p,q)}{d_{g_0}(p,q)} - d_{g_0}(p,q) \right) \\
                               & = \lim_{p \to \xi, q \to \eta} (d_{g_1}(p,q) - d_{g_0}(p,q)) \left( \frac{d_{g_1}(p,q) + d_{g_0}(p,q)}{d_{g_0}(p,q)} \right) \\
                               & = 2 S_{g_0}(g_1)(\xi, \eta) \\
\end{align*}

where in the last step we have used Lemmas \ref{schwarzlimit} and \ref{distbdd}. $\diamond$

\medskip

\section{The kernel of the ray transform}

\medskip

We keep the notation of the previous section. For $m \geq 1$ we denote by $\sigma : T^*X^{\otimes m} \to T^*X^{\odot m}$ the
symmetrization operator on $(0,m)$ tensors. Covariant differentiation with respect to the Levi-Civita connection of $g_0$ defines an operator
$\nabla^{g_0} : \Gamma(T^*X^{\otimes m}) \to \Gamma(T^*X^{\otimes m+1})$ (where for $E$ a vector bundle over $X$, $\Gamma(E)$ denotes
as usual the space of smooth sections of $E$). Restricting to symmetric tensors and composing with the
symmetrization operator, we obtain a differential operator acting on symmetric tensors,
$$
d^{g_0} := \sigma \circ \nabla^{g_0} : \Gamma(T^*X^{\odot m}) \to \Gamma(T^*X^{\odot m+1})
$$
For $m = 1$, if $u \in \Gamma(T^*X)$ is a smooth $1$-form, then the symmetric $(0,2)$-tensor $d^{g_0}v$ coincides
with the Lie derivative $\mathcal{L}_v g_0$ where $v \in \Gamma(TX)$ is the vector field dual to $u$ with respect
to the metric $g_0$. For $u \in \Gamma(T^*X^{\odot m})$ and a $g_0$-geodesic $\gamma : (a,b) \to X$, the following equality is valid:
$$
\frac{d}{dt} (u(\dot{\gamma}(t))) = 2 (d^{g_0}u)(\dot{\gamma}(t))
$$
%
%
%
%
%
It follows immediately that for a smooth compactly supported $1$-form $v$, we have $I_{g_0}(d^{g_0}v) = 0$.
Conversely we have the following:

\medskip

\begin{theorem} \label{kernel} Let $f \in \Gamma(T^*X^{\odot 2})$ be
a smooth symmetric $(0,2)$-tensor field with compact support.
If $I_{g_0}(f) = 0$ then $f = d^{g_0} v$ for some $v \in \Gamma(T^*X)$ with compact support. 
Moreover the support of $v$ only
depends on the support of $f$.
\end{theorem}

\medskip

The proof of the above theorem will follow easily from the
characterization of the kernel of the ray transform on a
completely dissipative Riemannian manifold given by Sharafutdinov
in \cite{sharafutdinov1}. In order to state his result,
we first recall the relevant notions from \cite{sharafutdinov1}.

\medskip

\begin{definition} A compact Riemannian manifold with boundary
$(M, g_0)$ is called a completely dissipative Riemannian manifold
(or CDRM for short) if the
following conditions are satisfied:

\smallskip

\noindent (1) The boundary $\partial M$ is strictly convex,
i.e. the second fundamental form of the boundary is positive-definite.

\smallskip

\noindent (2) For every $x \in M$ and every $\xi \in T_x M - \{0\}$
the maximal geodesic $\gamma$ with initial conditions
$\gamma(0) = x, \gamma'(0) = \xi$ is defined on a finite interval.
\end{definition}

\medskip

It is not hard to see in our situation that any
closed $g_0$-metric ball $M \subset X$ is a CDRM. We define for $M$ a CDRM the following:

$$
\partial T^1_{\pm}M := \{ (x,\xi) \in T^1 M | x \in \partial M, \pm<\xi, \nu(x)> \geq 0 \}
$$
where $\nu(x)$ denotes the outward normal to the boundary.

\medskip

\begin{definition} Let $(M, g_0)$ be a CDRM. The ray transform of $M$ is the linear operator $I_M$ defined by
\begin{align*}
I_M : \Gamma(T^*M^{\odot 2}) & \to C^{\infty}(\partial T^1_+ M) \\
                u            & \mapsto \left( I_M(u) : (x,\xi) \mapsto \int_{\tau_{-}(x,\xi)}^{0} u(\dot{\gamma}_{(x,\xi)}(s)) \ ds \right)
\end{align*}
where $\gamma_{(x,\xi)} : [\tau_{-}(x,\xi), 0] \to M$ is the maximal geodesic with initial conditions $\gamma_{(x,\xi)}(0) = x, \gamma_{(x,\xi)}'(0) = \xi$.
\end{definition}

\medskip

We may now state a version of Sharafutdinov's result suited to our purposes:

\medskip

\begin{theorem} Let $(M, g_0)$ be a CDRM of nonpositive sectional curvature and
let $f \in \Gamma(T^*M^{\odot 2})$. If $I_M(f) = 0$
then $f = d^{g_0}v$ for some $v \in \Gamma(T^*M)$ such that $v_{|\partial M} = 0$.
\end{theorem}

\medskip

\noindent{\bf Proof of Theorem \ref{kernel}:} Given $f \in \Gamma(T^*X^{\odot 2})$
with compact support such that $I_{g_0}(f) = 0$, choose
a closed $g_0$-metric ball $M$ containing the support of $f$
in its interior.
Then $M$ is a CDRM, and moreover, since the support of $f$
is contained in $M$, the integral of $f$ over any maximal geodesic of
$M$ coincides with the integral of $f$ over the
extension of the geodesic segment to a bi-infinite geodesic of $X$. Hence the equality
$I_{g_0}(f) = 0$ implies $I_M(f_{|M}) = 0$. Applying the previous theorem, there exists
$v \in \Gamma(T^*M)$ such that $f_{|M} = d^{g_0}v$ and $v_{|\partial M} = 0$.
Extending $v$ to be zero outside $M$ so that $v$ has compact support
the conclusion of the theorem follows. $\diamond$

\medskip

\section{Proofs of theorems}

\medskip

\subsection{Infinitesimal rigidity}

\medskip

\noindent{\bf Proof of Theorem \ref{isospectral}:} Since $S_{g_0}(g_t) = 0$ for $0 \leq t \leq 1$, it follows from
Theorem \ref{schwarzvarn} that $I_{g_t}(\dot{g_t}) = 0$ for $0 \leq t \leq 1$. By Theorem \ref{kernel}, this
implies existence of vector fields $v_t$ for $0 \leq t \leq 1$ such that $\dot{g_t} = \mathcal{L}_{v_t} g_t$. Moreover
by the hypothesis on the supports of $g_t - g_0$, it follows that the supports of the vector fields $v_t$ are contained
in a fixed compact. Hence we may integrate to obtain a $1$-parameter family of diffeomorphisms $f_t : X \to X$ which are equal
to the identity outside a fixed compact, such that $f^*_t g_t = g_0$ for $0 \leq t \leq 1$. $\diamond$

\medskip

\subsection{Local rigidity}

\medskip

We first recall some lemmas and notation from \cite{crokedairbekovsharafutdinov}.

\medskip

Let $(M, g_0)$ be a CDRM.
We denote by $\delta^{g_0}$ the {\it divergence operator} of the metric $g_0$ acting on
symmetric tensors, which is a first-order differential operator formally adjoint to the operator
$d^{g_0}$ (we refer to \cite{sharafutdinov1} for the precise definition). Symmetric tensors
$f$ such that $\delta^{g_0}(f) = 0$ are called "solenoidal".

\medskip

For $k \geq 1$ an integer and $0 < \alpha < 1$ a real number, we denote
by $C^{k,\alpha}(T^*M^{\odot 2})$ the space of $C^{k,\alpha}$-smooth symmetric $(0,2)$
tensor fields on $M$. Endowing $C^{k,\alpha}(T^*M^{\odot 2})$ with the natural $C^{k,\alpha}$
topology turns it into a topological Banach space, i.e. a topological vector space whose
topology can be defined by some norm making it a Banach space.

\medskip

We denote by Diff$^{k,\alpha}_0(M)$ the set of all $C^{k,\alpha}$-smooth
diffeomorphisms of M that are the identity on
the boundary, and endow Diff$^{k,\alpha}_0(M)$ with the natural $C^{k,\alpha}$ topology
(defined using some finite atlas; the resulting topology is independent of the choice of
atlas).

\medskip

Theorem 2.1 of \cite{crokedairbekovsharafutdinov} states:

\medskip

\begin{theorem}\label{shiftsol} \cite{crokedairbekovsharafutdinov}
For every neighborhood $U \subset $Diff$^{k,\alpha}_0(M)$ of the identity there
 is a neighborhood $W \subset C^{k,\alpha}(T^*M^{\odot 2})$ of
the metric tensor $g_0$ such that for every metric $g \in W$
there exists a diffeomorphism $\phi \in U$ for which the
tensor field ${\phi}^*g$ is solenoidal, i.e., $\delta^{g_0}({\phi}^*g) = 0$.
\end{theorem}

\medskip

The metric $g_0$ defines in the usual way inner products on all spaces of
tensor fields on $M$, in particular on $\Gamma(T^*M^{\odot 2})$; we denote
the inner product on this space by $(.,.)_{L^2(T^*M^{\odot 2})}$. Then Proposition
4.1 of \cite{crokedairbekovsharafutdinov} states:

\medskip

\begin{prop} \label{volestimate} \cite{crokedairbekovsharafutdinov} There is an $\epsilon > 0$ such that if
$f \in C^0(T^*M^{\odot 2})$ satisfies $||f||_{C^0(T^*M^{\odot 2})} < \epsilon$ and
$Vol_{g_0 + f}(M) \leq Vol_{g_0}(M)$ then
$$
(g_0,f)_{L^2(T^*M^{\odot 2})} \leq \frac{2}{3}||f||^2_{L^2(T^*M^{\odot 2})}
$$
\end{prop}

\medskip

The following lemma follows from the arguments in section 6 of \cite{crokedairbekovsharafutdinov}:

\medskip

\begin{lemma} \label{pestovestimate} \cite{crokedairbekovsharafutdinov} Suppose $M$
is nonpositively curved. There is a
constant $C > 0$ and a quadratic
first-order differential operator $L$ on functions
on $\partial T^1_{+} M$ such that for any $f \in C^{2,\alpha}(T^*M^{\odot 2})$ 
such that ${\delta}^{g_0}(f) = 0$ and $f_{|\partial M} = 0$,
$$
||f||^2_{L^2(T^*M^{\odot 2})} \leq C \int_{\partial T^1_{+} M} L(I_M(f)) d\Sigma^{2n-2}
$$
where $d\Sigma^{2n-2}$ denotes the natural Liouville volume form on the manifold
$\partial T^1_{+} M$ induced by the metric $g_0$.
\end{lemma}

\medskip

The arguments used to prove Lemma 6.1 of \cite{crokedairbekovsharafutdinov} yield the following:

\medskip

\begin{lemma}\label{postvest} \cite{crokedairbekovsharafutdinov} There is a constant $C > 0$
such that the following holds:

If $f \in C^{2,\alpha}(T^*M^{\odot 2})$ is such that $I_M(f) \geq 0$ on all of $\partial_+ T^1 M$,
then
$$
\int_{\partial T^1_{+} M} L(I_M(f)) d\Sigma^{2n-2} \leq C ||f||_{C^2} (g_0, f)_{L^2(T^*M^{\odot 2})}
$$
\end{lemma}

\medskip

We now have all the ingredients required to prove Theorem \ref{localrigidity}:

\medskip

\noindent{\bf Proof of Theorem \ref{localrigidity}:} Let $M$ be a closed $g_0$-ball
containing the given compact $K$ in its interior. Let $g$ be a metric on $X$ with sectional curvatures
bounded above by $-1$ such that the support of $g - g_0$ is contained in $K$, $\hat{id}_{g_0,g}$ is
Moebius and $Vol_{g}(K) = Vol_{g_0}(K)$. By Theorem
\ref{shiftsol}, we may choose an $\epsilon > 0$
such that if $||g - g_0||_{C^{2,\alpha}(T^*M^{\odot 2})} < \epsilon$ then
there is a diffeomorphism $\phi \in $Diff$^{2,\alpha}_0(M)$ such
that $\delta^{g_0}({\phi}^*g) = 0$ on $M$. Extending $\phi$ to be the identity outside
$M$, the same identity $\delta^{g_0}({\phi}^*g) = 0$ holds on $X$. Let $g_1$ be the metric
${\phi}^*g$ on $X$, then $\hat{id}_{g_0, g_1}$ is also Moebius (since $\phi = id$ outside $M$),
$Vol_{g_1}(M) = Vol_{g_0}(M)$, and the tensor $f := g_1 - g_0$ has support contained in $M$
and vanishes identically on $\partial M$.
It also follows from Theorem \ref{shiftsol} that given $\epsilon' > 0$ we can always choose
$\epsilon > 0$ small enough such that $f$ satisfies $||f||_{C^{2,\alpha}(T^*M^{\odot 2})} < \epsilon'$.
Thus choosing $\epsilon > 0$ small enough we may ensure that Proposition \ref{volestimate}
applies to $f$, so that the hypothesis $Vol_{g_1}(M) = Vol_{g_0}(M)$ implies that
$$
(g_0,f)_{L^2(T^*M^{\odot 2})} \leq \frac{2}{3}||f||^2_{L^2(T^*M^{\odot 2})}
$$

\medskip

Since the support of $f$ is contained in $M$, we have that
for any $v \in \partial T^1_{+} M$, if $\xi, \eta \in \partial_{g_0} X$ denote
the endpoints of the bi-infinite $g_0$-geodesic with initial velocity $v$, then
$$
I_M(f)(v) = I_{g_0}(f)(\xi,\eta) \geq 2 S_{g_0}(g_1)(\xi,\eta) = 0
$$
(using Lemma \ref{raybiggerschwarz} and the fact that $\hat{id}_{g_0, g_1}$ is Moebius).

\medskip

It now follows from Lemmas \ref{pestovestimate} and \ref{postvest} that there
are constants $C_1, C_2 > 0$ such that
\begin{align*}
||f||^2_{L^2(T^*M^{\odot 2})} & \leq C_1 \int_{\partial T^1_{+} M} L(I_M(f)) d\Sigma^{2n-2} \\
                              & \leq C_1 C_2 ||f||_{C^2(T^*M^{\odot 2})} (g_0, f)_{L^2(T^*M^{\odot 2})} \\
                              & \leq C_1 C_2 \epsilon \frac{2}{3} ||f||^2_{L^2(T^*M^{\odot 2})} \\
\end{align*}
thus choosing $\epsilon$ small enough so that $\epsilon \cdot C_1 C_2 \frac{2}{3} < 1$
implies that $||f||^2_{L^2(T^*M^{\odot 2})} = 0$, so $f = 0$ on $M$ and hence on all of $X$,
so $g_1 = g_0$, and $g$ is isometric to $g_0$. $\diamond$

\medskip

\bibliography{isomoeb}
\bibliographystyle{alpha}

\end{document}